\documentclass[letterpaper,12pt]{amsart}

\usepackage{amssymb}
\usepackage{graphicx}
\usepackage{pb-diagram,pb-xy} 
\usepackage{lscape}
\usepackage[cmtip,arrow]{xy} 
\usepackage[margin=1in]{geometry}  
\usepackage[dvips]{color} 

\numberwithin{equation}{section}

\theoremstyle{plain}
\newtheorem{theorem}[equation]{Theorem}

\newtheorem{proposition}[equation]{Proposition}
\newtheorem{prop}[equation]{Proposition}

\newtheorem*{claim*}{Claim}

\theoremstyle{definition}
\newtheorem{definition}[equation]{Definition}

\newtheorem{remark}[equation]{Remark}

\newtheorem{example}[equation]{Example}

\newcommand{\isom}{\cong}                       
\newcommand{\smsh}{\wedge}                      

\newcommand{\Smsh}{\bigwedge}                   

\newcommand{\cat}[1]{\mathcal{#1}}              

\DeclareMathOperator*{\colim}{colim}

\newcommand{\Hom}{\operatorname{Hom}}

\newcommand{\ord}[1]{$#1$\textsuperscript{th}}

\title[Composition product]{A note on the composition product of symmetric sequences}
\author{Michael Ching}
\address{Department of Mathematics \\
        Amherst College \\
        Amherst, MA 01002}
\email{mching@amherst.edu}

\begin{document}

\begin{abstract}
We consider the composition product of symmetric sequences in the case where the underlying symmetric monoidal structure does not commute with coproducts. Even though this composition product is not a monoidal structure on symmetric sequences, it has enough structure, namely that of a `normal oplax' monoidal product, to be able to define monoids (which are then operads on the underlying category) and make a bar construction. The main benefit of this work is in the dual setting, where it allows us to define a cobar construction for cooperads.
\end{abstract}

\maketitle
\thispagestyle{empty}

\section*{Introduction}

The category of symmetric sequences, in a closed symmetric monoidal category $(\cat{C},\smsh,S)$, is well known to have a monoidal product (called the composition product) whose monoids are precisely the operads in $\cat{C}$. This result is due to Smirnov \cite{smirnov:1982}. This perspective allows for standard techniques in the theory of monoids, in particular bar constructions, to be applied to operads. This result depends, however, on the assumption that the monoidal structure $\smsh$ on $\cat{C}$ is \emph{closed}, or more specifically, on the fact that $\smsh$ commutes with colimits. This condition is necessary in order that the composition product be associative.

On the other hand, operads can be defined in any symmetric monoidal category without the closed condition. It is therefore natural to wonder to what extent operads can still be treated as monoids when the composition product is not strictly associative. In particular, can we form a bar construction for operads in this wider context? In this paper, we show that the composition product is a `normal oplax monoidal product' in the sense of Day-Street \cite{day/street:2003} and that this slightly weaker structure allows, nonetheless, for sensible definitions of monoids and of bar constructions over them. In particular, it follows that operads in any symmetric monoidal category have a natural bar construction.

The main motivation for this paper is, in fact, \emph{co}operads. A cooperad in $\cat{C}$ can be considered as an operad in $\cat{C}^{op}$ (with the canonical symmetric monoidal structure on the opposite category). Although most interesting symmetric monoidal categories are closed symmetric monoidal (for example, based compactly-generated topological spaces, or any of the standard models for stable homotopy theory), their opposite categories are less often so. This is reflected in the fact that the monoidal product does not generally commute with finite \emph{products} in these topological examples. Our results therefore are necessary to define cobar constructions for cooperads in these categories.

Such cobar constructions were described previously by the author in \cite{ching:2005a}, in which the theory presented here was outlined vaguely. The present paper is intended to fill in the gaps in that presentation and to demonstrate the usefulness of considering normal oplax monoidal structures, and their monoids, in homotopy theory.

Here is an outline of the paper. In \S\ref{sec:pseudo} we describe what is meant by a normal oplax monoidal structure. In \S\ref{sec:comp} we show that the composition product forms part of a normal oplax monoidal structure on the category of symmetric sequences. In \S\ref{sec:monoids} we describe monoids for normal oplax monoidal structures, and note that an operad is precisely such an monoid for the composition product. Finally, in \S\ref{sec:bar} we describe the simplicial bar construction in this setting.

\section*{Acknowledgements}

This note is a clarification and extension of parts of my thesis \cite{ching:2005a}. Thanks to my advisor Haynes Miller for making me think more about this. Thanks also to various referees for bringing the notion of normal oplax monoidal structure to my attention and for other useful comments. This work was partially supported by NSF Award DMS-1144149.

\section{Normal oplax monoidal structures} \label{sec:pseudo}

In this section we define the notion of a normal oplax monoidal structure on a category. This structure appears in \cite{day/street:2003} as an example of a lax monoid (in the 2-cell dual of the 2-category of categories). It is a weakening of the notion of monoidal structure in the sense that the unitivity and associativity morphisms are not required to be isomorphisms. To deal with this change, we need explicit functors to stand for the higher iterates of the usual monoidal structure. These functors are then related by a set of associativity morphisms that satisfy an appropriate set of conditions analogous to the usual axioms for a monoidal structure.

\begin{definition} \label{def:oplax}
Let $\cat{E}$ be a category. A \emph{normal oplax monoidal structure} on $\cat{E}$ consists of the following data:
\begin{itemize}
  \item an object $I$ in $\cat{E}$ (which we refer to as the \emph{unit} of the structure);
  \item for each integer $n \geq 2$, a functor $\mu_n: \cat{E}^n \to \cat{E}$;
  \item for each triple of integers $(n,l,r)$ with $0 \leq l \leq l+r \leq n$, a natural transformation
  \[ \alpha_{n,l,r} : \mu_n(X_1, \dots, X_n) \to \mu_{n-r+1}(X_1, \dots, X_l, \mu_r(X_{l+1}, \dots, X_{l+r}), X_{l+r+1}, \dots, X_n). \]
\end{itemize}
To make sense of $\alpha_{n,l,r}$ when $r = 0,1,n$, we take $\mu_1: \cat{E} \to \cat{E}$ to be the identity functor, and $\mu_0: \cat{E}^0 \to \cat{E}$ to be the functor whose only value is the unit $I$. For example the map $\alpha_{n,l,0}$ takes the form
\[ \mu_n(X_1,\dots,X_n) \to \mu_{n+1}(X_1,\dots,X_l,I,X_{l+1},\dots,X_n). \]
The above data must satisfy the following conditions:
\begin{enumerate}
    \item for any values of $n$ and $l$, the maps $\alpha_{n,l,1}$ and $\alpha_{n,0,n}$ are the identity morphism on $\mu_n(X_1,\dots,X_n)$;
    \item for any values of $n,l,r,k,s$ that make sense, the diagram
    \[ \begin{diagram} \dgARROWLENGTH=3.5em
      \node{\scriptstyle \mu_n(X_1, \dots, X_n)} \arrow{e,t}{\alpha_{n,l,r}} \arrow{s,l}{\alpha_{n,k,s}} \node{\scriptstyle \mu_{n-r+1}(X_1, \dots, \mu_r(X_{l+1}, \dots, X_{l+r}), \dots, X_n)} \arrow{s,r}{\alpha_{n-r+1,k-r+1,s}} \\
      \node{\scriptstyle \mu_{n-s+1}(X_1, \dots, \mu_s(X_{k+1}, \dots, X_{k+s}), \dots, X_n)} \arrow{e,b}{\alpha_{n-s+1,l,r}} \node{\scriptstyle \mu_{n-r-s+2}(X_1, \dots, \mu_r(X_{l+1}, \dots, X_{l+r}), \dots, \mu_s(X_{k+1}, \dots, X_{k+s}), \dots, X_n)}
    \end{diagram} \]
    commutes;
    \item for any values of $n,l,r,k,s$ that make sense, the diagram
    \[ \begin{diagram} \dgARROWLENGTH=4em
       \node{\scriptstyle \mu_n(X_1, \dots, X_n)} \arrow{e,t}{\alpha_{n,l,r}} \arrow{s,l}{\alpha_{n,k,s}} \node{\scriptstyle \mu_{n-r+1}(X_1, \dots, \mu_r(X_{l+1}, \dots, X_{l+r}), \dots, X_n)} \arrow{s,l}{\mu_{n-r+1}(X_1, \dots, X_l, \alpha_{r,k-l,s}, X_{l+r+1}, \dots, X_n)} \\
      \node{\scriptstyle \mu_{n-s+1}(X_1, \dots, \mu_s(X_{k+1}, \dots, X_{k+s}), \dots, X_n)} \arrow{e,b}{\alpha_{n-s+1,l,r-s+1}} \node{\scriptstyle \mu_{n-r+1}(X_1, \dots, \mu_{r-s+1}(X_{l+1}, \dots, \mu_s(X_{k+1}, \dots, X_{k+s}), \dots, X_{l+r}), \dots, X_n)}
   \end{diagram} \]
   commutes.
\end{enumerate}
A \emph{normal oplax monoidal category} consists of a category $\cat{E}$ together with a normal oplax monoidal structure on $\cat{E}$.
\end{definition}

\begin{example} \label{ex:monoidal}
Let $\cat{E}$ be a category with a monoidal structure given by the bifunctor $\otimes: \cat{E}^2 \to \cat{E}$ with unit object $I$. Set
\[ \mu_n(X_1, \dots, X_n) := (\dots(X_1 \otimes X_2) \otimes \dots ) \otimes X_n \]
where $\mu_0$ takes the value $I$. The unit isomorphisms for the monoidal structure determine morphisms $\alpha_{n,l,0}$ and the associativity isomorphisms determine $\alpha_{n,l,r}$ for $r \geq 1$, as in Definition \ref{def:oplax}. The necessary conditions follow from the axioms for a monoidal structure. Thus a monoidal structure determines a normal oplax monoidal structure.

Conversely, suppose we have a normal oplax monoidal structure on $\cat{E}$ in which all the associativity morphisms $\alpha_{n,l,r}$ are isomorphisms. Then $\mu_2$ is a monoidal structure on $\cat{E}$ with associativity isomorphism
\[ \alpha_{3,1,2} \circ \alpha_{3,0,2}^{-1}: \mu_2(\mu_2(X,Y),Z) \to \mu_3(X,Y,Z) \to \mu_2(X,\mu_2(Y,Z)) \]
and unit isomorphisms
\[ \alpha_{1,0,0} : X \to \mu_2(I,X) \text{ and } \alpha_{1,1,0} : X \to \mu_2(X,I). \]

More generally, any non-symmetric cooperad $Q$ in a braided monoidal category $\cat{E}$ determines a normal oplax monoidal structure given by
\[ \mu_n(X_1,\dots,X_n) = Q(n) \otimes X_1 \otimes \dots \otimes X_n. \]
This is dual to Example 2.6 of \cite{batanin/weber:2011}.
\end{example}

\begin{remark} \label{rem:lax}
Day and Street \cite{day/street:2003} define a normal oplax monoidal structure with an associativity morphism of the form
\[ \beta_{(n_1,\dots,n_r)}: \mu_n(X_1,\dots,X_n) \to \mu_r(\mu_{n_1}(X_1,\dots,X_{r_1}),\dots,\mu_{n_r}(X_{r_1+\dots+r_{n-1}+1},\dots,X_n)) \]
for each ordered partition $n = n_1+\dots+n_r$. These determine our structure maps $\alpha_{n,l,r}$ by considering partitions in which all but one of the values $n_i$ is $1$ and recalling that $\mu_1$ is the identity functor. On the other hand we can reconstruct $\beta_{(n_1,\dots,n_r)}$ from the $\alpha_{n,l,r}$ as the composite
\[ \alpha_{r-1+n_r,r-1,n_r} \circ \alpha_{r-2+n_{r-1}+n_r,r-2,n_{r-1}} \circ \dots \circ \alpha_{1+n_2+\dots+n_r,1,n_2} \circ \alpha_{n,0,n_1}. \]
Via this correspondence, the conditions (1), (2) and (3) of Definition \ref{def:oplax} are equivalent to the unit and associativity conditions satisfied by the $\beta_{(n_1,\dots,n_r)}$.

In Day and Street's original context, the adjective `normal' refers to the restriction that $\mu_1$ is the identity, and we obtain an `oplax monoidal structure' by replacing this restriction with a natural transformation $\eta: 1_{\cat{E}} \to \mu_1$. In this case $\mu_1$ is a comonad on $\cat{E}$ via the maps $\beta_{(1)}$ and $\eta$.

There is also a dual notion of `lax monoidal structure' on a category $\cat{E}$. The only difference is that the directions of the associativity morphisms are reversed, that is, a (normal) lax monoidal structure on $\cat{E}$ is the same as a (normal) oplax monoidal structure on $\cat{E}^{op}$.
\end{remark}

\begin{remark} \label{rem:multicategory}
From a category $\cat{E}$ with normal oplax monoidal structure we can construct, in a natural way, a multicategory (see, for example, \cite[2.1.1]{leinster:2004}) with the same class of objects as $\cat{E}$. For $X_1,\dots,X_n,Y \in \cat{E}$, we define the set of multi-maps from $(X_1,\dots,X_n)$ to $Y$ by
\[ \Hom_{\cat{E}}(X_1,\dots,X_n;Y) := \Hom_{\cat{E}}(\mu_n(X_1,\dots,X_n),Y). \]
In this way, we can view a normal oplax monoidal structure as a special kind of multicategory: namely, one in which the functors
\[ Y \mapsto \Hom_{\cat{E}}(X_1,\dots,X_n;Y) \]
are represented, in a natural way, by single objects of $\cat{E}$.
\end{remark}

\section{Composition of symmetric sequences} \label{sec:comp}

Our main example of a normal oplax monoidal structure comes from the composition product on symmetric sequences.

\begin{definition} \label{def:symseq}
A \emph{symmetric sequence} in a category $\cat{C}$ is a functor $A : \mathsf{FinSet} \to \cat{C}$ from the category $\mathsf{FinSet}$, whose objects are finite sets and whose morphisms are bijections, to $\cat{C}$. Denote the category of all symmetric sequences in $\cat{C}$ by $\cat{C}^{\Sigma}$ (in which morphisms are natural transformations).

There is an equivalent way to define a symmetric sequence which explains the terminology. Let $\mathsf{FinSet}_0$ be the category whose objects are the finite sets $\underline{r} := \{1,\dots,r\}$ for $r \geq 0$ (with $\underline{0}$ the empty set), and whose morphisms are bijections. The category $\mathsf{FinSet}_0$ is a skeletal subcategory of $\mathsf{FinSet}$, so any symmetric sequence $A: \mathsf{FinSet} \to \cat{C}$ is determined, up to canonical isomorphism, by its restriction $A: \mathsf{FinSet}_0 \to \cat{C}$. This restriction consists of the sequence $A(\underline{0}),A(\underline{1}),A(\underline{2}),\dots$ of objects in $\cat{C}$, together with an action of the symmetric group $\Sigma_r$ on $A(\underline{r})$, hence the name `symmetric sequence'.
\end{definition}

\begin{remark}
For the remainder of this section, we fix a symmetric monoidal structure $\smsh$ with unit object $S$ on the category $\cat{C}$. The symmetry and associativity isomorphisms for the monoidal product $\smsh$ allow us to write expressions such as
\[ X \smsh Y \smsh Z \;\; \text{and} \;\; \bigwedge_{\alpha \in A} X_\alpha \]
without caring about parentheses or ordering of the factors.

We also assume throughout this paper that the underlying category $\cat{C}$ has all colimits (and, when we refer to the opposite category $\cat{C}^{op}$, we assume that $\cat{C}$ has all limits).
\end{remark}

\begin{remark}
In this paper we focus on the composition product for \emph{symmetric} sequences and our application is to \emph{symmetric} operads. It seems likely that there are corresponding results for non-symmetric sequences in an arbitrary braided (not necessarily symmetric) monoidal category $\cat{C}$. This would provide a simplicial bar construction for non-symmetric operads in the case that the monoidal structure on $\cat{C}$ is not closed.
\end{remark}

\begin{definition} \label{def:J/FinSet}
For a finite set $J$, we define a category $J/\mathsf{FinSet}_0$ as follows. The class of objects of $J/\mathsf{FinSet}_0$ consist of all functions $f:J \to I$ for some $I \in \mathsf{FinSet}_0$, and the set of morphisms from $(f:J \to I)$ to $(f':J \to I')$ is the set of bijections $\sigma: I \to I'$ such that $f' = \sigma \circ f$.
\end{definition}

\begin{definition}
Let $A$ and $B$ be two symmetric sequences in $\cat{C}$. For each finite set $J$, we define a functor
\[ (A,B): J/\mathsf{FinSet}_0 \to \cat{C} \]
on objects by
\[ (A,B)(f:J \to I) := A(I) \smsh \Smsh_{i \in I} B(f^{-1}(i)). \]
For a morphism $\sigma:I \to I'$ in $J/\mathsf{FinSet}_0$ we define
\[ (A,B)(\sigma) := A(I) \smsh \Smsh_{i \in I} B(f^{-1}(i)) \to A(I') \smsh \Smsh_{i' \in I'} B(f^{-1}(\sigma^{-1}(i'))) \]
by combining map $A(\sigma)$ with the permutation of the smash product identifying the term corresponding to $i \in I$ with the term corresponding to $\sigma(i) \in I'$.
\end{definition}

\begin{definition} \label{def:comp}
For symmetric sequences $A,B$, we define a symmetric sequence $A \circ B$ by
\[ (A \circ B)(J) := \colim_{f \in J/\mathsf{FinSet}_0} (A,B)(f). \]
A bijection $\theta: J \to J'$ determines a map $(A,B)(f) \to (A,B)(f \circ \theta^{-1})$ that is the map
\[ A(I) \smsh \Smsh_{i \in I} B(f^{-1}(i)) \to A(I) \smsh \Smsh_{i \in I} B(\theta(f^{-1}(i))) \]
via the identity on $A(I)$ and the action of the symmetric sequence $B$ on the bijections $f^{-1}(i) \isom \theta(f^{-1}(i))$ given by restricting $\theta$. We thus obtain induced maps
\[ (A \circ B)(\theta) := (A \circ B)(J) \to (A \circ B)(J') \]
that make $A \circ B$ into a symmetric sequence in $\cat{C}$.
\end{definition}

\begin{remark}
The definition of the composition product in Definition \ref{def:comp} is isomorphic to more familiar descriptions involving a coproduct over partitions of $J$.
Note that each element of $J/\mathsf{FinSet}_0$ amounts to a partition of $J$. (Note that here partitions are allowed to include multiple copies of the empty set. So even the empty set itself has infinitely many partitions, each corresponding to a function $\emptyset \to \underline{r}$.) Morphisms in $J/\mathsf{FinSet}_0$ connect only those partitions with the same number of pieces. The colimit $(A \circ B)(J)$ can therefore be written as
\[ (A \circ B)(J) \isom \coprod_{r = 0}^{\infty} \left[ \coprod_{J = J_1 \amalg \dots \amalg J_r} A(\underline{r}) \smsh B(J_1) \smsh \dots \smsh B(J_r) \right]_{\Sigma_r} \]
where $\Sigma_r$ acts by permuting the sets $J_1,\dots,J_r$ and on $A(\underline{r})$ in the usual way.
\end{remark}

\begin{definition} \label{def:unit}
The \emph{unit symmetric sequence} $\mathbb{I}$ in the symmetric monoidal category $(\cat{C},\smsh,S)$ is given by
\[ \mathbb{I}(J) = \begin{cases} S & \text{if $|J| = 1$}; \\ 0 & \text{otherwise}; \end{cases} \]
where $0$ is a fixed initial object in $\cat{C}$. The map $\mathbb{I}(J) \to \mathbb{I}(J')$ induced by a bijection $J \to J'$ is the identity morphism on $S$ or $0$ as appropriate.
\end{definition}

\begin{proposition} \label{prop:special}
Let $\cat{C}$ be a closed symmetric monoidal category. Then the composition product determines a monoidal structure on the category $\cat{C}^{\Sigma}$ of symmetric sequences in $\cat{C}$ with unit object $\mathbb{I}$.
\end{proposition}
\begin{proof}
See \S I.1.8 of \cite{markl/shnider/stasheff:2002}.
\end{proof}

The importance in Proposition \ref{prop:special} of the hypothesis that the monoidal structure be closed is that it implies that the monoidal product $\smsh$ has a right adjoint in each variable and hence commutes with colimits. This condition is needed to show that the composition product is associative. We now drop this assumption. Our goal then is to show that the composition product is now part of a normal oplax monoidal structure on $\cat{C}^{\Sigma}$. To do this we need to define the higher composition products of three or more symmetric sequences.

\begin{definition}
For a finite set $J$ and integer $n \geq 2$, we denote by $J/\mathsf{FinSet}_0[n]$ the category whose objects are sequences of functions of finite sets
\[ J \arrow{e,t}{f_{n-1}} I_{n-1} \arrow{e,t}{f_{n-2}} \dots  \arrow{e,t}{f_1} I_1 \]
with $I_1,\dots,I_{n-1} \in \mathsf{FinSet}_0$, and whose morphisms are commutative diagrams
\[ \begin{diagram}
  \node[2]{I_{n-1}} \arrow[2]{s,r}{\sigma_{n-1}} \arrow{e,t}{f_{n-2}} \node{I_{n-2}} \arrow{e,t}{f_{n-3}} \arrow[2]{s,r}{\sigma_{n-2}} \node{\dots} \arrow{e,t}{f_1} \node{I_1} \arrow[2]{s,r}{\sigma_{1}} \\
  \node{J} \arrow{ne,t}{f_{n-1}} \arrow{se,b}{f'_{n-1}} \\
  \node[2]{I'_{n-1}} \arrow{e,t}{f'_{n-2}} \node{I'_{n-2}} \arrow{e,t}{f'_{n-3}} \node{\dots} \arrow{e,t}{f'_1} \node{I'_1}
\end{diagram} \]
in which the vertical maps $\sigma_1,\dots,\sigma_{n-1}$ are bijections. Notice that $J/\mathsf{FinSet}_0[2]$ is the category $J/\mathsf{FinSet}_0$ of Definition \ref{def:J/FinSet}. We also define $J/\mathsf{FinSet}_0[1]$ to be the terminal category (with one object and one morphism).
\end{definition}

\begin{definition}
Now let $A_1, \dots, A_n$ be symmetric sequences in the symmetric monoidal category $\cat{C}$. We define a functor
\[ (A_1,\dots,A_n): J/\mathsf{FinSet}_0[n] \to \cat{C} \]
as follows. For the sequence
\[ \dgTEXTARROWLENGTH=2em f: J \arrow{e,t}{f_{n-1}} I_{n-1} \arrow{e,t}{f_{n-2}} \dots  \arrow{e,t}{f_1} I_1 \]
we set
\[ (A_1,\dots,A_n)(f) := A_1(I_1) \smsh \Smsh_{i \in I_1} A_2(f_1^{-1}(i)) \smsh \dots \smsh \Smsh_{i \in I_{n-1}} A_n(f_{n-1}^{-1}(i)). \]
For a morphism $\sigma: f \to f'$ in $J/\mathsf{FinSet}_0[n]$, the bijections $\sigma_k$ determine an isomorphism
\[ (A_1,\dots,A_n)(f) \to (A_1,\dots,A_n)(f'). \]
\end{definition}

\begin{definition} \label{def:highcomp}
Fix $n \geq 1$ and take symmetric sequences $A_1,\dots,A_n$ in $\cat{C}$. The \emph{higher composition product} of $A_1,\dots,A_n$ is the symmetric sequence $(A_1 \circ \dots \circ A_n)$ defined by
\[ (A_1 \circ \dots \circ A_n)(J) := \colim_{f \in J/\mathsf{FinSet}_0[n]} (A_1,\dots,A_n)(f). \]
For $n = 2$ this reduces to the ordinary composition product of Definition \ref{def:comp} and when $n = 1$ it produces the identity functor on the category of symmetric sequences. To define also the case $n = 0$, we define the composition product of the empty collection to be the unit symmetric sequence $\mathbb{I}$ of Definition \ref{def:unit}.
\end{definition}

\begin{remark}
As for the binary composition product, there is a formulation of the higher products in terms of partitions. Objects of the groupoid $J/\mathsf{FinSet}_0[n]$ are in one-to-one correspondence with sequences of nested partitions of $J$ of length $n-1$. If the pieces of the \ord{k} partition are labelled $J_{k,1}, \dots, J_{k,r_{k}}$ we have
\[ (A_1 \circ \dots \circ A_n)(J) \isom \coprod_{r_1,\dots,r_{n-1}} \left[ \coprod A_1(\underline{r_1}) \smsh \Smsh_{i = 1}^{r_1} A_2(J_{1,i}) \smsh \dots \smsh \Smsh_{i = 1}^{r_{n-1}} A_n(J_{n-1,i}) \right]_{\Sigma_{r_1} \times \dots \times \Sigma_{r_{n-1}}}. \]
\end{remark}

\begin{remark}
If $\smsh$ commutes with colimits in $\cat{C}$, the higher composition product $(A_1 \circ \dots \circ A_n)$ is isomorphic to any of the possible ways of iterating the binary product, for example,
\[ (A_1 \circ A_2 \circ A_3) \isom (A_1 \circ A_2) \circ A_3 \isom A_1 \circ (A_2 \circ A_3), \]
and this observation produces the associativity isomorphisms for the composition product in this case. Without this assumption, there are still maps from $(A_1 \circ A_2 \circ A_3)$ to each of these iterated composition products. We construct these next and then show that they form a normal oplax monoidal structure on $\cat{C}^{\Sigma}$.
\end{remark}

\begin{definition} \label{def:compassoc}
Let $A_1,\dots,A_n$ be symmetric sequences in the symmetric monoidal category $\cat{C}$. Our aim is to construct natural maps of symmetric sequences:
\[ \alpha_{n,l,r}: (A_1 \circ \dots \circ A_n) \to (A_1 \circ \dots \circ A_{l} \circ (A_{l+1} \circ \dots \circ A_{l+r}) \circ A_{l+r+1} \circ \dots \circ A_n). \]
For each finite set $J$, we need to give a map
\[ \colim_{f \in J/\mathsf{FinSet}_0[n]} (A_1,\dots,A_n)(f) \to \colim_{g \in J/\mathsf{FinSet}_0[n-r+1]} (A_1,\dots,A_l,(A_{l+1} \circ \dots \circ A_{l+r}),A_{l+r+1}, \dots, A_n)(g). \]
We do so by picking, for each $f \in J/\mathsf{FinSet}_0[n]$, an element $g \in J/\mathsf{FinSet}_0[n-r+1]$ and giving maps
\[ \alpha_{n,l,r}(f): (A_1,\dots,A_n)(f) \to (A_1,\dots,A_l,(A_{l+1} \circ \dots \circ A_{l+r}),A_{l+r+1}, \dots, A_n)(g) \]
that commute with the maps induced by a morphism $\sigma: f \to f'$ in $J/\mathsf{FinSet}_0[n]$.

So consider a sequence
\[ \dgTEXTARROWLENGTH=2em f: J \arrow{e,t}{f_{n-1}} I_{n-1} \arrow{e,t}{f_{n-2}} \dots  \arrow{e,t}{f_1} I_1 \]
and take $g \in \mathsf{FinSet}_0[n-r+1]$ to be sequence
\[ \dgTEXTARROWLENGTH=3em g: J \arrow{e,t}{f_{n-1}} \dots \arrow{e,t}{f_{l+r}} I_{l+r} \arrow{e,t}{f_{l+r-1} \dots f_{l}} I_{l} \arrow{e,t}{f_{l-1}} \dots  \arrow{e,t}{f_1} I_1. \]
Comparing the objects involved in the required map $\alpha_{n,l,r}(f)$, we see that it is sufficient to produce a map
\[ \left[ \left(\Smsh_{i \in I_l}A_{l+1}(f_l^{-1}(i))\right) \smsh \dots \smsh \left(\Smsh_{i \in I_{l+r-1}}A_{l+r}(f_{l+r-1}^{-1}(i)) \right) \right] \to \Smsh_{i \in I_l}(A_{l+1} \circ \dots \circ A_{l+r})((f_{l+r-1}\dots f_{l})^{-1}(i)). \]
But the left-hand side here is isomorphic to
\[ \Smsh_{i \in I_l} (A_{l+1},\dots,A_{l+r})(f[i]) \]
where $f[i]$ is the sequence
\[ \dgTEXTARROWLENGTH=3em (f_{l+r-1}\dots f_l)^{-1}(i) \arrow{e,t}{f_{l+r-1}} (f_{l+r-2}\dots f_l)^{-1}(i) \arrow{e,t}{f_{l+r-2}} \dots \arrow{e,t}{f_{l+1}} (f_l)^{-1}(i) \]
and so the required map $\alpha_{n,l,r}(f)$ is formed by taking the $\smsh$-product over $i \in I_l$ of the canonical maps
\[ (A_{l+1},\dots,A_{l+r})(f[i]) \to (A_{l+1} \circ \dots \circ A_{l+r})((f_{l+r-1}\dots f_l)^{-1}(i)). \]
These $\alpha_{n,l,r}(f)$ commute with the maps induced by a morphism $\sigma: f \to f'$ in $J/\mathsf{FinSet}_0[n]$ so define the natural transformation $\alpha_{n,l,r}$.
\end{definition}

\begin{theorem} \label{thm:comp}
Let $\cat{C}$ be a symmetric monoidal category with all colimits. Then the functors
\[ \mu_n(A_1,\dots,A_n) := (A_1 \circ \dots \circ A_n), \]
(with $\mu_1(A) = A$ and $\mu_0 = \mathbb{I}$, the unit symmetric sequence), and the maps $\alpha_{n,l,r}$ of Definition \ref{def:compassoc}, form a normal oplax monoidal structure on the category $\cat{C}^{\Sigma}$ of symmetric sequences in $\cat{C}$.
\end{theorem}
\begin{proof}
We check axioms (1)-(3) of Definition \ref{def:oplax}. For (1), we first look at Definition \ref{def:compassoc} for $\alpha_{n,l,1}$ and note that $g = f$ and $\alpha_{n,l,1}(f)$ is the identity map on $(A_1,\dots,A_n)(f)$. Thus $\alpha_{n,l,1}$ is the identity on $(A_1 \circ \dots \circ A_n)$. For $\alpha_{n,0,n}$ we notice that $g$ the unique object in $J/\mathsf{FinSet}_0[1]$, and that $\alpha_{n,l,r}(f)$ is the canonical map
\[ (A_1,\dots,A_n)(f) \to (A_1 \circ \dots \circ A_n)(J). \]
It follows that $\alpha_{n,0,n}$ is the identity.

For each of (2) and (3), we have to compare two maps originating in the colimit
\[ \colim_{f \in J/\mathsf{FinSet}_0[n]} (A_1,\dots,A_n)(f) \]
so it is sufficient to check that for each $f \in J/\mathsf{FinSet}_0[n]$, the components associated to $f$ by these two maps are equal. So consider a sequence
\[ \dgTEXTARROWLENGTH=2em f: J \arrow{e,t}{f_{n-1}} I_{n-1} \arrow{e,t}{f_{n-2}} \dots  \arrow{e,t}{f_1} I_1. \]
For (2) each of the relevant composites is built from maps
\[ (A_1,\dots,A_n)(f) \to (A_1,\dots,(A_{l+1} \circ \dots \circ A_{l+r}), \dots, (A_{k+1} \circ \dots \circ A_{k+s}), \dots, A_n)(h) \]
where $h$ is the sequence
\[ \dgTEXTARROWLENGTH=3em h: J \arrow{e,t}{f_{n-1}} \dots \arrow{e,t}{f_{k+s}} I_{k+s} \arrow{e,t}{f_{k+s-1} \dots f_{k}} I_{k} \arrow{e,t}{f_{k-1}} \dots \arrow{e,t}{f_{l+r}} I_{l+r} \arrow{e,t}{f_{l+r-1} \dots f_{l}} I_{l} \arrow{e,t}{f_{l-1}} \dots  \arrow{e,t}{f_1} I_1 \]
and each of these is built from the $\smsh$-product of canonical maps of the form
\[ (A_{l+1},\dots,A_{l+r})(f[i]) \to (A_{l+1} \circ \dots \circ A_{l+r})((f_{l+r-1}\dots f_l)^{-1}(i)) \]
for $i \in I_l$ and
\[ (A_{k+1},\dots,A_{k+s})(f[i']) \to (A_{k+1} \circ \dots \circ A_{k+s})((f_{k+s-1}\dots f_k)^{-1}(i')) \]
for $i' \in I_k$.

For (3) each composite is built from maps
\[ (A_1,\dots,A_n)(f) \to (A_1,\dots,(A_{l+1} \circ \dots \circ (A_{k+1} \circ \dots \circ A_{k+s}) \circ \dots \circ A_{l+r}),\dots,A_n)(g) \]
where $g$ is the sequence
\[ \dgTEXTARROWLENGTH=3em g: J \arrow{e,t}{f_{n-1}} \dots \arrow{e,t}{f_{l+r}} I_{l+r} \arrow{e,t}{f_{l+r-1} \dots f_{l}} I_{l} \arrow{e,t}{f_{l-1}} \dots  \arrow{e,t}{f_1} I_1. \]
Each of these is built from the $\smsh$-product of maps
\[ (A_{l+1},\dots,A_{l+r})(f[i]) \to (A_{l+1} \circ \dots \circ (A_{k+1} \circ \dots \circ A_{k+s}) \circ \dots \circ A_{l+r})(h[i]) \]
for $i \in I_l$ where $h[i]$ is the sequence
\[ \dgTEXTARROWLENGTH=3em (f_{l+r-1}\dots f_l)^{-1}(i) \arrow{e,t}{f_{l+r-1}} \dots (f_{k+s-1}\dots f_l)^{-1}(i) \arrow{e,t}{f_{k+s-1}\dots f_k} (f_{k-1} \dots f_l)^{-1}(i) \dots \arrow{e,t}{f_{l+1}} (f_l)^{-1}(i) \]
and each of these is built from the $\smsh$-product of the canonical maps
\[ (A_{k+1},\dots,A_{k+s})(h[i][i']) \to (A_{k+1} \circ \dots \circ A_{k+s})((f_{k+s-1} \dots f_k)^{-1}(i')) \]
for $i' \in I_k$.
\end{proof}

\section{Monoids in normal oplax monoidal categories} \label{sec:monoids}

We now return to the general context and say what is meant by a monoid in an arbitrary normal oplax monoidal structure on a category $\cat{E}$. We saw in Remark \ref{rem:multicategory} that such a structure determines a multicategory with the same class of objects. A monoid in $\cat{E}$ is then precisely a monoid in that multicategory in the sense of Leinster \cite[2.1.11]{leinster:2004}. This entails the following.

\begin{definition} \label{def:monoid}
Let $\cat{E}$ be a normal oplax monoidal category as in Definition \ref{def:oplax}. A \emph{monoid} with respect to this structure consists of an object $M \in \cat{E}$, morphisms
\[ m_n: \mu_n(M,\dots,M) \to M, \]
for $n \geq 0$, such that $m_1$ is the identity morphism on $M$, and such that for each triple $(n,l,r)$ with $0 \leq l \leq l+r \leq n$ the diagram
\[ \begin{diagram}
  \node{\mu_n(M,\dots,M)} \arrow{sse,b}{m_n} \arrow{e,t}{\alpha_{n,l,r}} \node{\mu_{n-r+1}(M,\dots,\mu_r(M,\dots,M),\dots,M)} \arrow{s,r}{m_r} \\
    \node[2]{\mu_{n-r+1}(M,\dots,M)} \arrow{s,r}{m_{n-r+1}} \\
    \node[2]{M}
\end{diagram} \]
commutes.
\end{definition}

\begin{example} \label{ex:monoid}
If the normal oplax monoidal structure on $\cat{E}$ comes from an actual monoidal structure as in Example \ref{ex:monoidal}, then a monoid in the sense of Definition \ref{def:monoid} is the same as a monoid for the monoidal structure in the usual sense.
\end{example}

\begin{remark} \label{rem:monoid-multicat}
A monoid in a multicategory $\cat{M}$ corresponds uniquely to a functor of multicategories $M: * \to \cat{M}$ where $*$ denotes the terminal multicategory, with a single object and a single multi-map of each possible type (see \cite[2.1.11]{leinster:2004}). A monoid $M$ in the normal oplax monoidal category $\cat{E}$ can similarly be described as a functor $* \to \cat{E}$ of normal oplax monoidal categories.
\end{remark}

We now show that a monoid $M$ in a normal oplax monoidal category is determined by the structure maps $m_2: \mu_2(M,M) \to M$ and $m_0: I \to M$. This is a consequence of the corresponding statement for monoids in a multicategory but we are unable to find such a result explicitly in the literature.

\begin{prop} \label{prop:monoid}
Let $\cat{E}$ be a normal oplax monoidal category an $M$ and object in $\cat{E}$. Suppose given maps
\[ m_2: \mu_2(M,M) \to M; \quad m_0: I \to M \]
such that the following diagrams commute
\begin{enumerate}
\item
\[ \begin{diagram}
  \node[2]{\mu_2(\mu_2(M,M),M)} \arrow{e,t}{\mu_2(m_2,M)} \node{\mu_2(M,M)} \arrow{se,t}{m_2} \\
  \node{\mu_3(M,M,M)} \arrow{ne,t}{\alpha_{3,0,2}} \arrow{se,b}{\alpha_{3,1,2}} \node[3]{M} \\
  \node[2]{\mu_2(M,\mu_2(M,M))} \arrow{e,t}{\mu_2(M,m_2)} \node{\mu_2(M,M)} \arrow{ne,b}{m_2}
\end{diagram} \]
\item
\[ \begin{diagram}
  \node{M} \arrow{e,t}{\alpha_{1,0,0}} \arrow{see,b}{1_M} \node{\mu_2(I,M)} \arrow{e,t}{\mu_2(m_0,M)} \node{\mu_2(M,M)} \arrow{s,r}{m_2} \\
  \node[3]{M}
\end{diagram} \]
\item
\[ \begin{diagram}
  \node{M} \arrow{e,t}{\alpha_{1,1,0}} \arrow{see,b}{1_M} \node{\mu_2(M,I)} \arrow{e,t}{\mu_2(M,m_0)} \node{\mu_2(M,M)} \arrow{s,r}{m_2} \\
  \node[3]{M}
\end{diagram} \]
\end{enumerate}
Then there is a unique monoid in $\cat{E}$ based on the object $M$ that includes the maps $m_2$ and $m_0$.
\end{prop}
\begin{proof}
The axioms for a monoid in Definition \ref{def:monoid} imply that for $n \geq 3$, $m_n$ is uniquely determined by $m_2$. This gives us the uniqueness part. For the existence claim, we take $m_1$ to be the identity on $M$, and define $m_3$ to be the composed map $\mu_3(M,M,M) \to M$ in diagram (1) above.

Our goal now is to show that it is possible to define $m_n$ for $n > 3$ in such a way that all diagrams of the form
\[ \tag{$D(n,l,r)$} \begin{diagram}
  \node{\mu_n(M,\dots,M)} \arrow{sse,b}{m_n} \arrow{e,t}{\alpha_{n,l,r}} \node{\mu_{n-r+1}(M,\dots,\mu_r(M,\dots,M),\dots,M)} \arrow{s,r}{m_r} \\
    \node[2]{\mu_{n-r+1}(M,\dots,M)} \arrow{s,r}{m_{n-r+1}} \\
    \node[2]{M}
\end{diagram} \]
commute. Notice that the cases $r = 1$ and $r = n$ are trivially satisfied (since $m_1$ is the identity). We focus first on the case $2 \leq r \leq n-1$ and save $r = 0$ for later.

Suppose that $n \geq 4$ and we have defined $m_p: \mu_p(M,\dots,M) \to M$ for $p < n$ such that diagrams $D(p,l,r)$, for $p < n$ and $2 \leq r \leq p-1$, commute. This is true for $n = 4$ by diagram (1) and our choice of $m_3$. For integers $l,r$ satisfying $2 \leq r \leq n-1$ and $r+l \leq n$, let
\[ \phi_{l,r}: \mu_n(M,\dots,M) \to M \]
be the composite appearing in diagram $D(n,l,r)$ (and to which we would like $m_n$ to be equal). We claim that
\[ \phi_{l,r} = \phi_{k,s} \]
for any two choices of $(l,r)$ and $(k,s)$. (We assume $l \leq k$.) To see this we write
\[ C_{l,r} := \{l+1,\dots,l+r\} \]
to refer to the sequence of copies of $M$ that are combined to form $\mu_r(M,\dots,M)$ in the diagram $D(n,l,r)$. The comparison between $\phi_{l,r}$ and $\phi_{k,s}$ depends on the relationship between the sequences $C_{l,r}$ and $C_{k,s}$.

When the sequences $C_{l,r}$ and $C_{k,s}$ are either disjoint or one is contained in the other, we consider the following diagram whose vertical composites are $\phi_{l,r}$ and $\phi_{k,s}$:
\[ \begin{diagram}
  \node[2]{\mu_n(M,\dots,M)} \arrow{sw,t}{\alpha_{n,l,r}} \arrow{se,t}{\alpha_{n,k,s}} \\
  \node{\mu_{n-r+1}(\dots,\mu_r(M,\dots,M),\dots)} \arrow{s,l}{m_r} \node[2]{\mu_{n-s+1}(\dots,\mu_s(M,\dots,M),\dots)} \arrow{s,r}{m_s} \\
  \node{\mu_{n-r+1}(M,\dots,M)} \arrow{se,b}{m_{n-r+1}} \node[2]{\mu_{n-s+1}(M,\dots,M)} \arrow{sw,b}{m_{n-s+1}} \\
  \node[2]{M}
\end{diagram} \]
This diagram can be filled in with either axiom (2) or (3) from Definition \ref{def:oplax} together with naturality squares and copies of diagrams $D(p,l,r)$ for $p < n$. Thus $\phi_{l,r} = \phi_{k,s}$.

Now suppose that the two sequences have intersection of length at least two or union of length at most $n-1$. Then, using the previous comparison of sequences for which one is contained in the other we have
\[ \phi_{l,r} = \phi_{k,l+r-k} = \phi_{k,s}. \]
Finally, suppose the intersection is length one and the union is length $n$. Since $n \geq 4$, we must have one of $r,s$ equal to at least $3$ in which case we either have
\[ \phi_{l,r} = \phi_{l,r-1} = \phi_{k,s} \]
or
\[ \phi_{l,r} = \phi_{k,s-1} = \phi_{k,s} \]
by the previous comparisons. The conclusion is that $\phi_{l,r} = \phi_{k,s}$ for all possibilities and so we define $m_n$ to be this common map. Recursively, this defines $m_n$ for all $n \geq 0$ such that all diagrams $D(n,l,r)$ for $r \geq 1$ commute.

It remains to show that with these choices each diagram $D(n,l,0)$ commutes. This takes the form
\[ \begin{diagram}
  \node{\mu_n(M,\dots,M)} \arrow{sse,b}{m_n} \arrow{e,t}{\alpha_{n,l,0}} \node{\mu_{n+1}(M,\dots,M,I,M,\dots,M)} \arrow{s,r}{m_0} \\
    \node[2]{\mu_{n+1}(M,\dots,M)} \arrow{s,r}{m_{n+1}} \\
    \node[2]{M}
\end{diagram} \]
When $l \neq n$, we can expand out the term $m_{n+1}$ using the diagram $D(n+1,l,2)$ to get
\[ \begin{diagram}
  \node{\mu_n(M,\dots, M, \dots,M)} \arrow{e,t}{\alpha_{1,0,0}} \arrow{sse} \node{\mu_n(M,\dots,\mu_2(I,M),\dots,M)} \arrow{s,r}{m_0} \\
    \node[2]{\mu_n(M,\dots,\mu_2(M,M),\dots,M)} \arrow{s,r}{m_2} \\
    \node[2]{\mu_n(M,\dots,M,\dots,M)} \arrow{s,r}{m_n} \\
    \node[2]{M}
\end{diagram} \]
The diagonal map is the identity on $\mu_n(M,\dots,M)$ by diagram (2) and so the overall composite is $m_n$ as required. Finally, if $l = n$, we use a similar argument based on diagram (3) instead.
\end{proof}

\begin{example} \label{ex:operad}
Let $(\cat{C},\smsh,S)$ be a symmetric monoidal category. A monoid in the normal oplax monoidal category $\cat{C}^{\Sigma}$ of Theorem \ref{thm:comp} is precisely an operad in $\cat{C}$. This generalizes the result of Smirnov \cite{smirnov:1982} to the case that the symmetric monoidal structure on $\cat{C}$ is not closed.
\end{example}

We can also talk about actions of a monoid in a normal oplax monoidal category. Again these are examples of corresponding notions for multicategories but these notions do not seem to be explicitly described in the literature (though are no doubt well known). We spell out the details.

\begin{definition} \label{def:module}
Let $M$ be a monoid in the normal oplax monoidal category $\cat{E}$ with structure maps $m_n$. A \emph{left $M$-module} consists of an object $L \in \cat{E}$ and, for each $n \geq 1$ a morphism
\[ l_n : \mu_n(M,\dots,M,L) \to L \]
(where $l_1: L \to L$ is the identity morphism) such that each of the following diagrams commutes
\begin{enumerate}
\item for $l+r < n$
\[ \begin{diagram}
  \node{\mu_n(M,\dots,M,L)} \arrow{sse,b}{l_n} \arrow{e,t}{\alpha_{n,l,r}} \node{\mu_{n-r+1}(M,\dots,\mu_r(M,\dots,M),\dots,M,L)} \arrow{s,r}{m_r} \\
  \node[2]{\mu_{n-r+1}(M,\dots,M,L)} \arrow{s,r}{l_{n-r+1}} \\
  \node[2]{L}
\end{diagram} \]
\item for $l + r = n$
\[ \begin{diagram}
  \node{\mu_n(M,\dots,M,L)} \arrow{sse,b}{l_n} \arrow{e,t}{\alpha_{n,l,r}} \node{\mu_{n-r+1}(M,\dots,M,\mu_r(M,\dots,M,L))} \arrow{s,r}{l_r} \\
  \node[2]{\mu_{n-r+1}(M,\dots,M,L)} \arrow{s,r}{l_{n-r+1}} \\
  \node[2]{L}
\end{diagram} \]
\end{enumerate}
The corresponding notion of \emph{right $M$-module} consists of an object $R \in \cat{E}$ with
\[ r_n : \mu_n(R,M,\dots,M) \to R \]
satisfying similar conditions.
\end{definition}

\begin{remark}
The structure of a left (or right) $M$-module is uniquely determined by the map $l_2: \mu_2(M,L) \to L$ (respectively, $r_2$) subject to conditions similar to diagrams (1) and (2) (respectively, (1) and (3)) of Proposition \ref{prop:monoid} with the final copy of $M$ replaced by $L$ (respectively, the first copy replaced by $R$).

The combination of a monoid $M$ and left $M$-module $L$ in $\cat{E}$ can be described as a map of multicategories
\[ F: (m \blacktriangleright l) \to \cat{E} \]
where $(m \blacktriangleright l)$ is the multicategory with two objects and exactly one multi-map of each of the types
\[ (m,\dots,m) \to m \]
and
\[ (m,\dots,m,l) \to l. \]
There is a similar characterization of right modules.
\end{remark}

\begin{example}
If $M$ is a monoid in the normal oplax monoidal category $\cat{C}^{\Sigma}$ of Theorem \ref{thm:comp} (that is, $M$ is an operad in $\cat{C}$), then a left or right $M$-module is a left or right module over the operad $M$ in the usual sense of, for example, \cite[3.26]{markl/shnider/stasheff:2002}.
\end{example}

\section{Simplicial bar construction on a monoid} \label{sec:bar}

In this section we show that we can form a simplicial bar construction for a monoid and its modules in any normal oplax monoidal category. This generalizes the usual construction for a monoid in a monoidal category.

\begin{definition} \label{def:bar}
Let $M$ be a monoid with respect to a normal oplax monoidal structure on a category $\cat{E}$. Let $L$ be a left $M$-module and $R$ a right $M$-module. For $n \geq 0$, we define
\[ \mathcal{B}_n(R,M,L) := \mu_{n+2}(R, \underbrace{M, \dots, M}_n, L) \]
Now let $\theta: [n] \to [m]$ be an order-preserving map (where $[n]$ denotes the totally-ordered set $\{0,\dots,n\}$). Then we define
\[ \theta^*: \mathcal{B}_m(R,M,L) \to \mathcal{B}_n(R,M,L) \]
to be the composite
\[ \begin{diagram}
    \node{\mu_{m+2}(R,M,\dots,M,L)} \arrow{e,t}{\alpha} \arrow{se,b}{\theta^*} \node{\mu_{n+2}(\mu_{\theta_0}(R,\dots),\dots,\mu_{\theta_j}(M,\dots,M),\dots,\mu_{\theta_{n+1}}(\dots,L))} \arrow{s,r}{(r_{\theta_0},\dots,l_{\theta_{n+1}})} \\
    \node[2]{\mu_{n+2}(R,M,\dots,M,L)}
\end{diagram} \]
where
\[ \theta_i := \begin{cases}
    \quad \theta(0) & \text{if $i = 0$}; \\
    \theta(i) - \theta(i-1) & \text{if $1 \leq i \leq n$}; \\
    m - \theta(n) & \text{if $i = n+1$}.
\end{cases} \]
The map marked $\alpha$ is a composite of the associativity maps of the normal oplax monoidal structure on $\cat{E}$, and that marked $(r_{\theta_0},\dots,l_{\theta_{n+1}})$ is built from the monoid and module action maps for $R$, $M$ and $L$.
\end{definition}

\begin{theorem} \label{thm:simpbar}
Let $\cat{E}$ be a normal oplax monoidal category with monoid $M$, left $M$-module $L$ and right $M$-module $R$. Definition \ref{def:bar} then determines a simplicial object $\mathcal{B}_{\bullet}(R,M,L)$ in $\cat{E}$ which we refer to as the \emph{simplicial bar construction} for $R,M,L$.
\end{theorem}
\begin{proof}
We must show that for $\theta: [n] \to [m]$ and $\zeta: [m] \to [l]$ we have
\[ (\zeta \theta)^* = \theta^* \zeta^*. \]
This amounts to the commutativity of the following diagram (notation condensed to fit):
\[ \begin{diagram} \dgARROWLENGTH=1em
  \node{\scriptstyle \mu_{l+2}(R,\dots,L)} \arrow{e,t}{\alpha} \arrow{s,l}{\alpha} \node{\scriptstyle \mu_{n+2}(\mu_{(\zeta \theta)_0}(R,\dots),\dots,\mu_{(\zeta \theta)_{n+1}}(\dots,L))} \arrow{se,t}{(r_{(\zeta \theta)_0},\dots,l_{(\zeta \theta)_{n+1}})} \arrow{s,l}{\alpha} \\
  \node{\scriptstyle \mu_{m+2}(\mu_{\zeta_0}(R,\dots,),\dots,\dots,\mu_{\zeta_{m+1}}(\dots,L))} \arrow{e,t}{\alpha} \arrow{s,l}{(r_{\zeta_0},\dots,l_{\zeta_{m+1}})} \node{\scriptstyle \mu_{n+2}(\mu_{\theta_0}(\mu_{\zeta_0}(R,\dots),\dots),\dots,\mu_{\theta_{n+1}}(\dots,\mu_{\zeta_{m+1}}(\dots,L)))} \arrow{s,l}{(r_{\zeta_0},\dots,l_{\zeta_{m+1}})} \node{\scriptstyle \mu_{n+2}(R,\dots,L)} \\
  \node{\scriptstyle \mu_{m+2}(R,\dots,L)} \arrow{e,b}{\alpha} \node{\scriptstyle \mu_{n+2}(\mu_{\theta_0}(R,\dots),\dots,\mu_{\theta_{n+1}}(\dots,L))} \arrow{ne,b}{(r_{\theta_0},\dots,l_{\theta_{n+1}})}
\end{diagram} \]
in which maps marked $\alpha$ are composites of the associativity maps for the normal oplax monoidal structure, and those marked $(r,\dots,l)$ involve the module and monoid structure maps for $L$, $R$ and $M$. The top-left square commutes by combining copies of the axioms of Definition \ref{def:oplax}, the bottom-left is a naturality square, and the quadrilateral on the right commutes by combining copies of the axioms of Definitions \ref{def:module} and \ref{def:monoid}.
\end{proof}

\begin{example}
When the normal oplax monoidal structure on $\cat{E}$ arises from an actual monoidal product $\otimes$ on $\cat{E}$ as in Example \ref{ex:monoidal}, the bar construction $\mathcal{B}_\bullet(R,M,L)$ is given by
\[ \mathcal{B}_n(R,M,L) = R \otimes \underbrace{M \otimes \dots \otimes M}_n \otimes L \]
with the usual simplicial structure.
\end{example}

\begin{example}
Now let $\cat{C}^\Sigma$ be the category of symmetric sequences on a symmetric monoidal category $\cat{C}$, with the normal oplax monoidal structure given by Theorem \ref{thm:comp}. Then the simplicial bar construction described in this section is the `standard' simplicial bar construction on an operad $P$ with respect to a left $P$-module $L$ and right $P$-module $R$. See, for example \cite[7.9]{ching:2005a}.

If we replace $\cat{C}$ with its opposite category $\cat{C}^{op}$, together with the canonical symmetric monoidal structure, the simplicial bar construction on an operad in $\cat{C}^{op}$ gives us a \emph{cosimplicial cobar construction} $\mathfrak{C}^\bullet(R,Q,L)$ on a cooperad $Q$ in $\cat{C}$ with respect to a right $Q$-comodule $R$ and left $Q$-comodule $L$. Explicitly:
\[ \mathfrak{C}^n(R,Q,L)(J) := \lim_{f \in J/\mathsf{FinSet}_0[n]^{op}} (R,\underbrace{Q,\dots,Q}_n, L)(f). \]
\end{example}

\bibliographystyle{amsplain}
\bibliography{C:/Users/mching/Dropbox/WinEdt/mcching}

\providecommand{\bysame}{\leavevmode\hbox to3em{\hrulefill}\thinspace}
\providecommand{\MR}{\relax\ifhmode\unskip\space\fi MR }
\providecommand{\MRhref}[2]{%
  \href{http://www.ams.org/mathscinet-getitem?mr=#1}{#2}
}
\providecommand{\href}[2]{#2}
\begin{thebibliography}{1}

\bibitem{batanin/weber:2011}
Michael Batanin and Mark Weber, \emph{Algebras of higher operads as enriched
  categories}, Applied Categorical Structures \textbf{19} (2011), no.~1,
  93--135.

\bibitem{ching:2005a}
Michael Ching, \emph{Bar constructions for topological operads and the
  {G}oodwillie derivatives of the identity}, Geom. Topol. \textbf{9} (2005),
  833--933 (electronic). \MR{MR2140994}

\bibitem{day/street:2003}
Brian Day and Ross Street, \emph{Lax monoids, pseudo-operads, and convolution},
  Diagrammatic morphisms and applications ({S}an {F}rancisco, {CA}, 2000),
  Contemp. Math., vol. 318, Amer. Math. Soc., Providence, RI, 2003, pp.~75--96.
  \MR{1973511 (2004c:18011)}

\bibitem{leinster:2004}
Tom Leinster, \emph{Higher operads, higher categories}, London Mathematical
  Society Lecture Note Series, vol. 298, Cambridge University Press, Cambridge,
  2004. \MR{2094071 (2005h:18030)}

\bibitem{markl/shnider/stasheff:2002}
Martin Markl, Steve Shnider, and Jim Stasheff, \emph{Operads in algebra,
  topology and physics}, Mathematical Surveys and Monographs, vol.~96, American
  Mathematical Society, Providence, RI, 2002. \MR{2003f:18011}

\bibitem{smirnov:1982}
V.~A. Smirnov, \emph{On the cochain complex of topological spaces}, Math. USSR
  Sbornik \textbf{43} (1982), 133--144.

\end{thebibliography}

\end{document}